\documentstyle[12pt]{article}
\newcommand{\natu}{\mbox{I$\!$N}}
\newcommand{\real}{\mbox{I$\!$R}}

\newcommand{\ext}{e^{\theta (r)}}
\newcommand{\G}{G (\theta (r))\,}
\newcommand{\dext}{(e^{\theta (r)})'}
\newcommand{\qed}{{\hfill {$\rlap{$\sqcap$}\sqcup$}}\\[0.2in]\hspace*{0.5in}}
\newcommand{\bk}{\\[0.03in] \hspace*{0.5in} }
\begin{document}

\title{Rotationally Symmetric $F$-Harmonic Maps Equations}
\author{Man Chun Leung\\Department of Mathematics, National University of
Singapore,\\ Singapore 119260  \ \ {\tt {matlmc@nus.sg}}}
\date{May, 1996}
\maketitle
\begin{abstract}We study a second order differential equation corresponding
to rotationally symmetric $F$-harmonic maps. We show unique continuation and
Liouville's type theorems for positive solutions. Asymptotic properties 
and the
existence of bounded positive solutions are investigated.
\end{abstract}

\vspace{0.5in}

{\bf \Large 1. \ \ Introduction}

\vspace{0.3in}

In this paper we study positive solutions $\alpha$ of the differential 
equation
\begin{eqnarray*}
(1.1) \ \ \ \ \ \ \ \G \alpha'' (r) & + & \left[  (n - 1) \G  {{f' 
(r)}\over {f
(r)}} + {d\over {dr}} \G \right] \alpha' (r) \\ & \  & \ \ \ - (n - 1) \G {{g
(\alpha (r)) g' (\alpha (r))}\over {f^2 (r)}} = 0\,, \ \ \ \ \ \ \ r > 
0\,, \ \
\ \ \ \
\end{eqnarray*}
with
$$\lim_{r \to 0^+} \alpha (r) = 0\,,$$
where  
$$
\theta (r) = {1\over 2} \left[ (\alpha' (r))^2 + (n - 1) {{g^2 (\alpha 
(r))}\over
{f^2 (r)}} \right]\,. \leqno (1.2)
$$
Equation (1.1) arises as the Euler-Lagrange equation of a
functional on the energy of rotationally symmetric maps. Let
$n
\ge 2$ be an integer and $(S^{n - 1}, d\vartheta^2)$ be the unit sphere in
${\real}^n$ with the induced Riemannian metric $d\vartheta^2$. Consider the
following {\it model} Riemannian manifolds \cite{G-W} 
$$
 M (f) = ( [0, \infty ) \times S^{n-1}\,, \ dr^2 + f^2(r) d\vartheta^2 
\ )\,,
$$      
$$  N (g) = ( [0, \infty ) \times S^{n-1}\,, \ dr^2 + g^2(r)d\vartheta^2\ 
)\,,
$$   
where $f, \,g \in C^2 ([0, \infty))$ satisfy the following conditions 
$$  f (0) = g (0) = 0\,, \ \ \ \ f' (0) = g' (0) = 1\,, \ \ \ 
f (r) > 0 \ \ {\mbox{and}} \ \ g(r) > 0
\leqno (1.3)
$$
for all $r > 0\,.$ $M (f)$ and $N (g)$ are complete noncompact Riemannian
manifolds. In particular, the Euclidean space and the hyperbolic space are
corresponding to
$f (r) = r$ and $g (r) = \sinh r$, respectively. A map $\Phi_\alpha : M 
(f) \to N
(g)$ is called a rotationally symmetric map if 
$$
\Phi_\alpha \,(r, \vartheta) = (\alpha (r)\,, \vartheta) \ \ \ \ 
{\mbox{for \ \
all}}
\
\
\ r > 0 \ \ \ {\mbox{and}} \ \ \vartheta \in S^{n - 1}\,,
$$   
where $\alpha : [0,
\infty) \to [0, \infty)$ is a function with $\alpha (0) = 0$.  Assume that
$\alpha
\in C^2 ((0, \infty))$. The energy density \cite{Eells-L-1} of a rotationally
symmetric map
$\Phi_\alpha : M (f)
\to M (g)$ is given by 
$$
\theta (r) = {1\over 2} \left[ (\alpha'
(r))^2 + (n - 1) {{g^2 (\alpha (r))}\over {f^2 (r)}} \right]
$$
for $r > 0$. Let $F : [0, \infty) \to [0, \infty)$ be a
function in $C^1 ([0, \infty)) \cap C^3 ((0, \infty))\,.$  Assume that
$F$ satisfies the following conditions
$$
F (x) > 0\,, \ \ {{dF}\over {dx}} (x) > 0 \ \ \ {\mbox {and}} \ \ \ {{d^2 
F}\over
{dx^2}} (x) 
\ge 0
\
\
\
\ {\mbox{for \ \ all}} \ \ x > 0\,. \leqno (1.4)
$$
Given a rotationally symmetric map $\Phi_\alpha : M (f) \to N (g)$, 
consider the
functional 
$$
\int_0^\infty F (\theta (r)) f^{n - 1} (r) \,dr\,. \leqno (1.5)
$$
Let 
$$G (x) = {{d F}\over {dx}} (x)\,, \ \ \ \ \ x \ge 0\,.$$
A standard variational argument shows that the Euler-Lagrange equation 
for the
functional (1.5) is given by equation (1.1). A rotationally symmetric map
$\Phi_\alpha$ is said to be a {\it rotationally symmetric $F$-harmonic 
map} if it
is a critical of the functional in (1.5), or equivalently, if $\alpha$ 
satisfies
equation (1.1) with
$\lim_{r \to 0^+}
\alpha (r) = 0\,.$ There are some interesting special
cases.\\[0.1in] 
({\bf A}) \ \ If $F (x) = x$
for all $x \ge 0$, then equation (1.1) becomes
$$\alpha''(r) + (n - 1) {{f'(r)}\over {f(r)}} \alpha' (r)  -  (n - 1){{ g
(\alpha (r)) g' (\alpha (r))}\over {f^2 (r)}} = 0\,, \leqno (1.6)$$
which is the equation for rotationally symmetric harmonic maps between model
Riemannian manifolds $M (f)$ and $N (g)$. We refer to the comprehensive
papers by Eells and Lemaire \cite{Eells-L-1} \cite{Eells-L-2} for 
harmonic maps.
Ratto and Rigoli \cite{R-R} study the asymptotic properties of positive
solutions to equation (1.6). Cheung and Law \cite{Ch-Law} use an initial 
value
approach to study equation (1.6).\\[0.1in]  
({\bf B})
\
\ If
$F (x) = x^{p \over 2}$ for some number
$p > 2$ and for all
$x \ge 0$, then equation (1.1) becomes
\begin{eqnarray*} (1.7) \ \ \ \ \  \ \ \theta^{{p\over 2} - 1} (r) 
\alpha''(r) &
+ & [(n - 1)
\theta^{{p\over 2} - 1}(r) {{f'(r)}\over {f(r)}} +  (\theta^{{p\over 2} -
1} (r))' ] \alpha'(r)\\ & \ & \ \ \ \ \ \ \ \ \ \ \ \ \ \ -  (n - 1)
\theta^{{p\over 2} - 1} (r){{ g (\alpha (r)) g' (\alpha (r))}\over {f^2 
(r)}} =
0\,,\ \ \ \ \ \ \ \ \ \ \ \ \  \ \ 
\end{eqnarray*}
which is the equation for rotationally symmetric $p$-harmonic maps 
between  $M
(f)$ and $N (g)$. Equation (1.7) is studied in \cite{Ch-Law-Leung}
\cite{L1} \cite{L2}.\\[0.1in] ({\bf C}) \ \ If $F (x) = e^x$ for all $x \ge
0$, then equation (1.1) becomes
\begin{eqnarray*} (1.8) \ \ \ \ \  \ \ \ext \alpha''(r) & + & [(n - 1)
\ext {{f'(r)}\over {f(r)}} +  \dext ] \alpha'(r)\\ & \ & \ \ \ \ \ \ \ \ 
\ \ \ \
\ \ -  (n - 1)
\ext {{ g (\alpha (r)) g' (\alpha (r))}\over {f^2 (r)}} = 0\,,\ \ \ \ \ \ 
\ \ \ \
\
\ \  \ \ 
\end{eqnarray*}
which is the equation for rotationally symmetric exponentially harmonic maps
between  $M (f)$ and $N (g)$. Exponentially harmonic maps are studied by
Eells and Lemaire \cite{Eells-L-3}, Duc and Eells \cite{Duc-Eells} and Hong
\cite{Hong}. We refer to \cite{L3} for properties of positive solutions to
equation (1.8).\bk  
In this paper we study local and
asymptotic properties for positive solutions to equation (1.1), where
$F$ is a $C^2$-function which satisfies the conditions in (1.4). We show
that a positive solution to equation (1.1) cannot go down to zero too 
fast as $r
\to 0^+$. This unique continuation property is well-known for harmonic maps
(cf. \cite{Eells-L-1}), but is unknown for $p$-harmonic maps (cf.
\cite{H-K-M}) and for general
$F$-harmonic maps. We show a uniqueness theorem for solutions of equation 
(1.1) 
under the assumption that $g'' (y) \ge 0$ for all $y > 0$, that is, the 
radial
Ricci curvature of the Riemannian manifold
$N (g)$ is nonpositive. Nonexistence of bounded positive
solutions to equation (1.1) is discussed. Roughly speaking, if $M (f)$ is 
"small"
and
$N (g)$ is comparatively large, then any bounded rotationally symmetric
$F$-harmonic map has to be the constant map. We find that under growth
conditions on
$F$, the derivative of solutions to equation (1.1) can be bounded.\bk 
In section 4, we discuss
the existence of bounded rotationally symmetric
$F$-harmonic maps from $M (f)$ to $N (g)$. Our results show that any
rotationally symmetric $F$-harmonic maps from the hyperbolic space to the
Euclidean space have to be bounded. In case $M (f)$ and $N (g)$ are the
hyperbolic space, we show that if $\alpha (r_o) < r_o$ for some point 
$r_o >
0\,,$ then $\alpha$ is a bounded function.\bk  
Our study provides a uniform approach to some properties for the 
equations in
(A), (B) and (C), although each equation has its own distinctive features.
Equation (1.1) is a quasilinear second order differential equation. We 
introduce
a first order equation on the energy function $\theta (r)$ (equation (2.6)),
which we can use to obtain useful information.

\vspace{0.5in}

{\bf \Large 2. \ \ Local Properties}

\vspace{0.3in}

In this section we assume that $\alpha \in C^2 ((0, R))$ is a
solution to equation (1.1) with
$\lim_{r \to 0^+} \alpha (r) = 0\,,$ where $R$ is a positive number,
and $F : [0, \infty) \to [0, \infty)$ satisfies the conditions in
(1.4).\\[0.2in]  {\bf Lemma 2.1.}
\
\ {\it For
$n
\ge 2$, assume that
$\alpha (r) > 0$ for all $r
\in (0, R)$ and
$g' (y) > 0$ for all $y \in {\real}^+,$ then $\alpha' (r) > 0$ for all $r \in
(0, R)\,.$}\\[0.1in]
{\bf Proof.} \ \ As
$\alpha (r) > 0$ for all $r \in (0, R)$ and $\lim_{r \to 0^+} \alpha (r) 
= 0$, we
can find a point
$r_o \in (0, R)$ such that $\alpha ' (r_o) > 0$. Suppose that there 
exists a
point
$r_1 > r_o$ such that $\alpha' (r_1) \le 0$, then we can find a point $r' \in
(r_o, r_1]$ such that 
$$
\alpha (r') > 0\,, \ \ \alpha ' (r') = 0 \ \ {\mbox{and}} \ \ \alpha '' (r')
\le  0\,.
$$ 
At the point 
$r'$, equation (1.1) shows that
$$ 
G (\theta (r')) \alpha''(r') = (n - 1) G (\theta (r')) {{ g (\alpha (r')) g'
(\alpha (r'))}\over {f^2 (r')}}\,.
$$  
Since $g' (\alpha (r)) > 0$ for all $r \in (0, R)\,,$ and $\theta (r') > 0$,
therefore $G (\theta (r')) > 0$. We have $\alpha'' (r') > 0$, contradiction.
Therefore
$\alpha' (r) > 0$ for all
$r
\in (r_o\,, R)$. As $\lim_{r \to 0^+} \alpha (r) = 0\,,$ we can let $r_o
\to 0$. {\hfill {$\rlap{$\sqcap$}\sqcup$}}\\[0.2in]
{\bf Lemma 2.2.} \ \ {\it For $n = 1$, assume that $\alpha (r) \ge 0$ for 
all $r
\in (0, R)$. Then $\alpha (r) = c r$ for all $r \in (0, R)\,,$ where $c$ 
is a
non-negative constant.}\\[0.1in]
{\bf Proof.} \ \  Equation (1.1) gives
$$\left( \G \alpha' (r) \right)' = 0\,,$$
that is, there exists a constant $c_o$ such that 
$$
\G \alpha' (r) = c_o
$$
for all $r \in (0, R)$. If $\alpha \not\equiv 0$ on $(0, R)\,,$ then we can
find a point $r'$ with $\alpha (r') > 0$ and $\alpha' (r') > 0\,.$ 
Therefore $G
(\theta (r')) > 0$ and hence $c_o$ is a positive constant. We have 
$$ \G  = {{c_o}\over {\alpha' (r) }} \leqno (2.3)$$
for all $r \in (0, R)$. As in this case 
$\theta (r) = (\alpha' (r))^2/2$, we have
$$
{{d \G}\over {d r}} = {{d G }\over {d \theta}} ({\theta (r)}) \,\alpha' (r)
\alpha'' (r) \leqno (2.4)
$$
and by (2.3) we have
$$
{{d \G}\over {d r}} = -{{c_o}\over {(\alpha' (r))^2}} \alpha'' (r)\,. \leqno
(2.5)
$$
Therefore if $\alpha'' (r) \not= 0$, then (2.4) and (2.5) imply that 
$${{d G }\over {d \theta}} (\theta (r)) < 0$$
which contradicts the assumption in (1.4). Thus $\alpha'' (r) = 0$ for 
all $r
\in (0, R)\,,$ that is, $\alpha (r) = c r$ for some constant $c \ge
0\,.${\hfill {$\rlap{$\sqcap$}\sqcup$}}\\[0.2in]\hspace*{0.5in}
Differentiating equation (1.2) we have
\begin{eqnarray*}
\G {{d \theta (r)}\over {dr}} & = & \G [ \alpha' (r)  \alpha'' (r) + 
(n - 1)  {{g (\alpha (r) ) g' (\alpha (r))
\alpha' (r)}\over {f^2 (r)}}\\
& \ & \ \ \  -  (n - 1) {{g^2 (\alpha (r)) f' (r)}\over {f^3
(r)}} ]\,.
\end{eqnarray*}
Using the equation $d \G /dr = (d G/d \theta) \theta'$ and equation (1.1) to
cancel the term
$G \alpha''$ we obtain
\begin{eqnarray*} 
(2.6) \ & \ &  \theta' (r) \left[ \G + {{d G}\over {dx}}
\vert_{x = \theta (r)} (\alpha' (r))^2 \right] \\ & = & (n - 1) \G \{ {{2 g
(\alpha (r) ) g' (\alpha (r))
\alpha' (r)}\over {f^2 (r)}} - {{f' (r)}\over {f (r)}} [ {{ g^2 (\alpha
(r))}\over {f^2 (r)}} + (\alpha' (r))^2 ] \}\,.  \ \ \ \ \ \ \  \ \ \ \ 
\end{eqnarray*}
\vspace{0.05in}

{\bf Lemma 2.7.} \ \ {\it For $n \ge 1$, assume that $\alpha \ge 0$ on $(0,
R)$. If $\alpha (r) = O (r^k)$ near
$0$ for some $k > (n - 1)/2 + 1$, then $\alpha \equiv 0$ on
$(0, R)\,.$}\\[0.1in] {\bf Proof.} \ \ If $n  = 1$, then by lemma 2.2 we have
$\alpha (r) = c r$ for some nonnegative constant $c$. If $\alpha (r) = O 
(r^k)$
near
$0$ for some $k > 1\,,$ then $c = 0$ and hence $\alpha \equiv 0$ on $(0, R)$.
Assume that $n \ge 2$. Given
$\epsilon
\in (0, R)$ small enough and
$\kappa
\in \natu$, by (1.3) we can find positive constants
$a_o, b_o$ and $c_o$ depending on $\epsilon$ and $\kappa$ such that
$$  |f' (r) | \le a_o\,, \ \ \ f (r) \ge c_or \ \ \ {\mbox{and}} \ \ \ 0 
< g'
(\alpha (r)) < b_o \ \ \ \ {\mbox{for \ \ all}} \ \ \ r \in (0, {{\kappa\over
{\kappa + 1}}
\epsilon})\,.
$$   
From (2.6) we obtain
\begin{eqnarray*} (2.8)  & \ & \ \theta' (r) \left[ \G + {{d
G}\over {dx}} \vert_{x = \theta (r)} (\alpha' (r))^2 \right] \\  & \le & 
(n - 1)
\G
\left\{ {{b_o}\over {f (r)}} [ {{ g^2 (\alpha (r))}\over {f^2 (r)}} + 
(\alpha'
(r))^2 ] + {{a_o}\over {f (r)}}  [ {{ g^2 (\alpha (r))}\over {f^2 (r)}} +
(\alpha' (r))^2 ] \right\}\\ & \le & (n - 1) \G {{(a_o + b_o)}\over 
{c_or}}[ {{
g^2 (\alpha (r))}\over {f^2 (r)}} + (\alpha' (r))^2 ] \ \ \ {\mbox{for \ \
all}} \ \ r \in (0, {{\kappa\over {\kappa + 1}}
\epsilon})\,.
\end{eqnarray*} 
Assume that $\alpha \not\equiv 0$ on
$(0, \epsilon /2)$. We claim that $\alpha$ cannot be zero on $(0, 
\delta)$ for
any $\delta \in (0,
\epsilon /2)$. Suppose that $\alpha \equiv 0$ on
$(0, \delta)$ for some $\delta \in (0, \epsilon /2)$. Since $\alpha 
\not\equiv 0$
on $(0, \epsilon /2 )$, we may assume that $\alpha (r) > 0$ on $(\delta,
\epsilon /2)$. Thus $\G > 0$ on $(\delta\,, \epsilon /2)\,.$ From (2.8) 
we have 
$$
\theta' (r) \le {C\over r} \theta (r) \ \ \ \ \ \ {\mbox{for \ \ all}} \ 
\ \ r
\in (0, {{\kappa\over {\kappa + 1}} \epsilon})\,,\leqno (2.9)
$$  
where $C = C (a_o, b_o, c_o, \epsilon\,, n)$ is a positive
constant. Integrating (2.9) we have
$$\ln \theta (r) |^a_{b} \le C \ln r |^a_b\,,$$ where $\delta < b < a < 
\epsilon
/2$. That is 
$$ 
\theta (a) \le \theta (b) ({a\over b})^C\,.\leqno (2.10)
$$ 
Let $b \to \delta > 0$, we have $\theta (b) \to 0$, but $\theta (a) > 0$,
which contradicts (2.10). Therefore 
$\alpha$ cannot be zero on $(0, \delta)$ for any
$\delta \in (0, \epsilon /2)$. Thus there exists a point
$r_o \in (0, \epsilon /2)$ such that $\alpha (r_o) > 0$ and $\alpha' 
(r_o) > 0$.
Suppose that there exists a point $r_1 \in (0, \epsilon /2)$
$r_1 > r_o$ such that $\alpha (r_1) = 0$, then there exists a point $r' \in
(r_o, r_1)$ such that $\alpha (r') > 0$, $\alpha' (r') = 0$ and $\alpha'' (r)
\le 0$. At $r'$, equation (1.1) gives
$$ \G \alpha''(r') = (n - 1) \G {{ g (\alpha (r')) g' (\alpha (r'))}\over 
{f^2
(r)}}\,.$$ Therefore we have $\alpha'' (r') > 0$, contradiction. Hence 
$\alpha (r) > 0$ for all $r \in (r_o, \epsilon /2)$. As $\lim_{r \to 0^+} 
\alpha
(r) = 0$ and
$\alpha
 \not\equiv 0$ on $(0, \delta)$ for all $\delta > 0$, we can let
$r_o \to 0$. Thus $\alpha (r) > 0$ on $(0,
\epsilon /2)$. As in lemma 2.1 we can show that $\alpha' (r) > 0$ on $(0,
\epsilon /2)$. Given any $\delta > 0$, since $\lim_{r \to 0^+} \alpha (r) 
= 0$,
we can find
$\epsilon_o < \epsilon$ such that 
$$ 0 < \alpha (r) < \delta \ \ \ \ {\mbox {on}} \ \ \ \ (0, \epsilon_o)\,.
$$  As
$f, g \in C^2 ([0, \infty))$ with $f' (0) = g' (0) = 1$, given
$\varepsilon_1 \in (0, 1)$, we can assume that
$\delta$ and $\epsilon_o$ small enough such that 
$$ 
0 < f' (r)  \le 1 + \varepsilon_1\,, \ \ \ \  f (r) \ge (1 - 
\varepsilon_1) r
\ \ \ {\mbox{and}} \ \ \ 0 < g' (\alpha (r)) < 1 + \varepsilon_1\ \ 
{\mbox{for \ \ all}} \ r \in (0, \epsilon_o)\,.
$$ 
As in (2.9), using (2.6) we have
$$
\theta' (r) \le (n - 1) ({{1 + \varepsilon_1}\over {1 -
\varepsilon_1}} ) \theta (r) \ \ \ \ {\mbox {on}} \ \ \ \ (0, \epsilon_o)\,.
\leqno (2.11)
$$   
Let 
$$ C = (n - 1) {{(1 + \varepsilon_1)}\over {(1 - \varepsilon_1)}}\,.
$$ 
Integrating (2.11) we obtain 
$$ 
\theta (a) \le \theta (r) ({a\over r})^C\,,
$$  where in this case $\epsilon_o  > a > r > 0$ and $a$ is a constant. 
We have
$$
\theta (r) \ge C' r^C 
$$   
for some positive constant $C'$.  That is,
$$
\theta (r) = (n - 1) {{g^2 (\alpha (r))}\over {f^2 (r)}} + (\alpha' 
(r))^2 \ge
C' r^C\,.\leqno (2.12)
$$   
We may assume that $\epsilon_1$ is so small such that  
$$ 
 k > {C \over 2} + 1\,.
$$    
Since $\alpha (r) \le C_k r^k$ for all $r \in (0, \epsilon_o)\,,$ where
$C_k$ is a positive constant, and $f (r), g(r) \sim r$ when $r$ is small, 
we have
$$ (n - 1) {{g^2 (\alpha (r))}\over {f^2 (r)}} \le C'' r^{2k - 2}\,,
$$ where $C''$ is a positive constant.  Hence by (2.11) we have
$$ (\alpha' (r))^2 \ge C_o r^C \leqno (2.13)
$$  
for all $r \in (0, \epsilon_o)$. Here $C_o$ is a positive constant. 
Therefore we
have
$$
\alpha (r) - \alpha (b) = \int^r_b \alpha' (s) ds \ge \sqrt{C_o} \int_b^r
s^{{C}\over {2}} ds = {{2\sqrt{C_o} }\over {C + 1}} (r^{ {{C}\over {2}} + 
1 } -
b^{ {{C}\over {2}} + 1})\,, \ \ \
\epsilon_o > r > b > 0\,.
$$  As $\lim_{b \to 0^+} \alpha (b) = 0$, we have
$$
\alpha (r) \ge {{2\sqrt{C_o} }\over {C + 1}} r^{ {{C}\over {2}} +1 } \ \ 
\ \
{\mbox{for \ \ all}} \ \ \ r \in (0, \epsilon_o)\,,
\leqno (2.14)
$$ 
contradicting that $\alpha (r) \sim O (r^k)$ for some $k > C/2 + 1$.\qed
It follows from the proof of lemma 2.1 and lemma 2.7 that if $\alpha \ge 
0$ and
$\alpha \not\equiv 0$ on $(0, R)\,,$ then $\alpha (r) > 0$ and $\alpha' 
(r) >
0$ for all $r \in (0, R)$. We have the following uniqueness
theorem.\\[0.2in]  {\bf Theorem 2.15.} \ \ {\it Let
$\alpha\,,
\beta 
\in C^2 (0, R)$ be positive solutions to equation (1.1) for some $R > 0$ and
$\lim_{r
\to 0^+}
\alpha (r) = \lim_{r \to 0^+} \beta (r) = 0$. Assume that $g'' (y) \ge 0$ 
for all
$y \in {\real}^+$. If there is a point $R_o \in (0, R)$ such that $\alpha 
(R_o)
= \beta (R_o)$, then $\alpha \equiv
\beta$ on $(0, R)$.}\\[0.1in]  {\bf Proof.} \ \ For $n = 1$, we have 
$\alpha (r)
= c_o r$ and $\beta (r) = c_1 r$ for all $r \in (0, R)$, where $c_o$ and 
$c_1$
are positive numbers. If there is a point $R_o \in (0, R)$ such that $\alpha
(R_o) = \beta (R_o)$, then
$\alpha \equiv \beta$ on $(0, R)$. We may assume that $n \ge 2$. Suppose that
$\alpha \not\equiv
\beta$ on $(0, R_o)$. We may assume that there is a point $r_1 < R_o$ 
such that
$\alpha (r_1) > \beta (r_1)$. Let $r_o$ be the biggest number in $[0, 
r_1)$ such
that 
$$\lim_{r \to r_o^+} \alpha (r) = \lim_{r \to r_o^+} \beta (r)\,.$$ If 
$\alpha'
(r) \le \beta' (r)$ for all $r \in (r_o, r_1)$, then an integration shows 
that
$\alpha (r_1) \le \beta (r_1)$. Therefore we can find a point $r'
\in (r_o, r_1)$ such that 
$$
\alpha (r') > \beta (r') \ \ \ \ {\mbox{and}} \ \ \ \ \alpha' (r') >\beta'
(r')\,. \leqno (2.16)
$$   As $g'' (y) \ge 0$ for all $y > 0$ and $g' (0) = 1$, $g$ and
$g'$ are positive and nondecreasing functions. Lemma 2.1 implies that 
$\alpha'
(r) > 0$ and $\beta' (r) > 0$ for all $r \in (0, R)$. For $r \in (0, R)$ 
let 
$$\theta_\alpha (r) =  {1\over 2} \left[ \alpha'^2 (r) + (n - 1) {{g^2 
(\alpha
(r))}\over {f^2 (r)}} \right]\,,$$
$$\theta_\beta (r) =  {1\over 2} \left[ \beta'^2 (r) + (n - 1) {{g^2 (\beta
(r))}\over {f^2 (r)}} \right]\,.$$ As $G$ is nondecreasing, (2.16) implies
that
$$
\ln \, [ G (\theta_\alpha (r'))\, \alpha' (r') ] > \ln \, [ G (\theta_\beta
(r'))\,
 \beta' (r') ]\,.
\leqno (2.17)
$$ 
Let $R' \in (r', R_o]$ be the smallest number such that
$\alpha (R') =
\beta (R')$. Therefore $\alpha (r) > \beta (r)$ for all $r \in (r', R')$. A
similar argument as above shows that
$\alpha' (R') \le \beta' (R')$. In this case $R' > 0$ and $f$ is positive 
in a
neighborhood of
$R'$. If $\alpha' (R') = \beta' (R')$ and $\alpha (R') =
\beta (R')$, then the uniqueness result for equation (1.1) (c.f. 
\cite{C}, p.
259) implies that
$\alpha \equiv \beta$ on a neighborhood of $R'$. A continuation argument 
shows
that
$\alpha
\equiv
\beta$ on
$(0, R)$. Therefore we may assume that
$\alpha' (R') <
\beta' (R').$ Hence we have
$$
\ln \, [ G (\theta_\alpha (R'))\, \alpha' (R')] < \ln \, [ G (\theta_\beta
(R')) \,
 \beta' (R') ]\,.
\leqno (2.18)
$$ 
The inequality (2.17) and (2.18) imply that  there exists a point
$\bar r \in (r'\,, R']$ such that 
$$
\ln \, [ G (\theta_\alpha (\bar r)) \alpha' (\bar r) ]   <
 \ln \, [ G (\theta_\beta (\bar r))
 \beta' (\bar r) ] \leqno (2.19)
$$ 
and
$$
\left( \ln \, [ G (\theta_\alpha (\bar r)) \alpha' (\bar r) ] \right)' 
\le 
\left(  \ln \, [ G (\theta_\beta (\bar r))
 \beta' (\bar r) ] \right)'\,, \leqno (2.20)
$$ 
since $\bar r \in (r'\,, R']$, we have $\alpha (\bar r) \ge \beta (\bar 
r)\,.$
From equation (2.19) we have
$$
G (\theta_\alpha (\bar r)) \alpha' (\bar r) < G (\theta_\beta (\bar r)) 
\beta'
(\bar r)\,. \leqno (2.21)
$$ 
As $G$ is nondecreasing, we have $\alpha' (\bar r) < \beta (\bar r)\,.$
Using equation (1.1) we obtain
\begin{eqnarray*}
(2.22) \ \ \ \ \ \ \ \ \ \left( \ln \, [ G (\theta_\alpha (\bar r)) 
\alpha' (\bar
r) ]
\right)' & = & {{G (\theta (\bar r)) \alpha'' (\bar r) + {d\over {dr}} G 
(\theta
(\bar r))
\alpha' (\bar r)}\over {G (\theta (\bar r)) \alpha' (\bar r)}}\\
 & = & (n - 1)\left\{ {{ g (\alpha (\bar r)) g' (\alpha (\bar r))}\over {f^2
(\bar r)
\alpha' (\bar r)}} - {{f' (\bar r)}\over {f (\bar r)}} \right\}\,. 
 \ \ \ \ \ \ \ \ \ \ \  \ \ \ \ \ \ \ \ \ \ \ 
\end{eqnarray*}  
Similarly
$$
\left( \ln \, [ G (\theta_\beta (\bar r)) \beta' (\bar r) ] \right)' =   
(n -
1)\left\{ {{ g (\beta (\bar r)) g' (\beta (\bar r))}\over {f^2 (r)
\beta' (\bar r)}} - {{f' (\bar r)}\over {f (\bar r)}} \right\}\,. 
 \leqno (2.23)
$$ 
We have $\alpha (\bar r) \ge  \beta (\bar r)$ and $0 < \alpha' (\bar r) < 
\beta'
(\bar r)\,.$ As $g$ and $g'$ are nondecreasing positive functions, (2.22) and
(2.23) gives
$$
\left( \ln \, [ G (\theta_\alpha (\bar r)) \alpha' (\bar r) ] \right)' > 
\left(
\ln \, [ G (\theta_\beta (\bar r)) \beta' (\bar r) ]  \right)'\,,
$$ which contradicts (2.20). Therefore $\alpha' (R') = \beta' (R')$ and
hence $\alpha
\equiv
\beta$ on
$(0, R)\,.${\hfill {$\rlap{$\sqcap$}\sqcup$}}\\[0.2in]\hspace*{0.5in}
As the radial Ricci curvature of the Riemannian manifold $N (g)$ is given by
$- (n - 1) g'' (y)/ g (y)$ for $n \ge 2$ and $y > 0$, the assumption that 
$g''
(y) \ge 0$ for all $y > 0$ is equivalent to the condition that the radial
Ricci curvature of $N (g)$ is nonpositive. We note that uniqueness 
results for 
harmonic maps often require the curvature to be nonpositive (cf.
\cite{Ha} \cite{J-K}).\\[0.2in] 
{\bf Proposition 2.24.} \ \ {\it Assume that in equation (1.1)
$f \equiv g$ on
$[0, \infty)$ and $f'' (y) = g'' (y) \ge 0$ for all $y \in {\real}^+$. If 
$\alpha
\in C^2 ((0\,,
R))$ is a positive solution to equation (1.1) with $\lim_{r \to 0^+} \alpha
(r) = 0$, then either $\alpha (r) < r$ or $\alpha > r$ or $\alpha (r) = 
r$ for
all $r \in (0, R)$.}\\[0.1in] {\bf Proof.} \ \ As the function $\alpha_1 (r)
= r$ for all $r \in {\real}^+$ is a positive solution to equation (1.1) 
with $f
\equiv g\,.$ The result follows from theorem 2.15.{\hfill
{$\rlap{$\sqcap$}\sqcup$}}

\vspace{0.5in}

{\bf \Large 3. \ \ Liouville's type theorems}

\vspace{0.3in}

Let $\alpha \in C^2 ((0, \infty))$ be a negative solution to equation (1.1)
with $\lim_{r \to 0^+} \alpha (r) = 0\,.$ We have the following
Liouville's type theorem, under the assumption that $G$ grows up
polynomially or at most exponentially.\\[0.2in] 
{\bf Theorem 3.1.}
\
\ {\it Assume that there exist positive numbers $b$ and $c$ such that 
$|f' (r)|
\le b\,,$
$g' (r) > 0$ and
$g'^2 (r) - g (r) g'' (r)
\ge c^2$ for all $r > 0$. Assume either} 
$$G (x) \ge c \,x \,{{d G}\over {dx}} (x) \ \ \ \ {\mbox{for \ \ all}} \ 
\ x > 0
\leqno (A)$$
or 
$$
G (x) \ge c_o \,{{d G}\over {dx}} (x) \ \ \ \ {\mbox{for \ \ all}} \ \ x 
> 0\,,
\leqno (B)
$$
{\it where $c$ and $c_o$ are positive constants. If $\alpha$ is bounded, then
$\alpha \equiv 0$ on $(0, \infty)$.}\\[0.1in] 
{\bf Proof.} \ \ By lemma 2.2 we may assume that $n \ge 2$. Suppose that 
$\alpha
\not\equiv 0\,.$ It follows from lemma 2.1 and the proof of lemma 2.7 that
$\alpha (r) > 0$ and
$\alpha' (r) > 0$ 
for all
$r
\in {\real}^+$. As
$\alpha$ is bounded, we can find a positive constant $a$ such that $|g' 
(\alpha
(r))| \le a$ for all $r \in {\real}^+$. We have
\begin{eqnarray*} (3.2) \ \ \ \ \ \ \ \ \ \ & \ & {{2 g (\alpha (r) ) g' 
(\alpha
(r))
\alpha' (r)}\over {f^2 (r)}} - {{f' (r)}\over {f (r)}} [ {{ g^2 (\alpha
(r))}\over {f^2 (r)}} + (\alpha' (r))^2 ] \\ & \le & {1\over {f (r)}} \{ a\,
[{{g^2 (\alpha (r) )}\over {f^2 (r)}} + (\alpha' (r))^2] + b\, [{{g^2 
(\alpha (r)
)}\over {f^2 (r)}} + (\alpha' (r))^2] \} \\  & \le & {1\over {f (r)}} (a 
+ b) 
[{{g^2 (\alpha (r) )}\over {f^2 (r)}} + (\alpha' (r))^2]\,.  \ \ \ \ \ \ 
\ \ \ \ 
\ \ \ \ \ \ \ \ \ \ \ \ \ \ \ \ \ \ \ \ \ \ \ \ \ \ \ \ \ \ \ 
\end{eqnarray*} 
By (2.6) we have
\begin{eqnarray*}
(3.3) \ \ \ \ \ \ \ \ \ \ \ \ \  \ \ \ & \  & \theta' (r) \left( \G + {{d 
G}\over
{dx}}|_{x =
\theta (r)} (\alpha' (r))^2 \right)\\ 
& \le & (n - 1) {{\G}\over {f (r)}} (a + b)  [{{g^2 (\alpha (r) )}\over {f^2
(r)}} + (\alpha' (r))^2]\,. \ \ \ \ \ \ \ \ \ \ \ \ \ \ \ \ \ \ \ \ \ \ 
\end{eqnarray*}
If we assume that (A) $G (x) \ge c x {{dG}\over {dx}} (x)$ for all $x > 0$,
then we have
\begin{eqnarray*} 
& \  & \theta' (r) {{d G}\over {dx}}|_{x =
\theta (r)}  \left[ c \theta (r)  + (\alpha' (r))^2 \right]\\  & \le & (n
- 1) {{\G}\over {f (r)}} (a + b)  [{{g^2 (\alpha (r) )}\over {f^2 (r)}} +
(\alpha' (r))^2]\,.
\end{eqnarray*}
Therefore we can find a positive constant $C$ such that 
$$
{d\over {dr}} \G = \theta' (r) {{d G}\over {dx}}|_{x =
\theta (r)} < C {{\G}\over {f (r)}}\,. \leqno (3.4)
$$
If we assume that (B) $G (x) \ge c_o{{dG}\over {dx}} (x)$ for all $x > 0$,
then we have
\begin{eqnarray*}  & \  & \theta' (r) {{d G}\over {dx}}|_{x =
\theta (r)}  ( c_o   + (\alpha' (r))^2 )\\  & \le & (n - 1) {{\G}\over
{f (r)}} (a + b)  [{{g^2 (\alpha (r) )}\over {f^2 (r)}} + (\alpha' (r))^2]\,.
\end{eqnarray*}

As $\alpha$ is a bounded function and $f (r) \ge c_1^2$ for all $r \ge 1$,
where $c_1$ is a positive constant, we can find a positive constant $C_1$ 
such
that  
$$ {{[{{g^2 (\alpha (r) )}\over {f^2 (r)}} + (\alpha' (r))^2]}\over {[ 
c_o +
(\alpha' (r))^2] }} \le C_1\
$$
for all $r \ge 1\,.$
Thus we also have 
$$ {d\over {dr}} \G = \theta' (r) {{d G}\over {dx}}|_{x =
\theta (r)} < C {{\G}\over {f (r)}}\leqno (3.5)
$$
for all $r \ge 1$ and for some constant $C > 0$. Replace the term 
${d\over {dr}}
\G$ in equation (1.1) by the above inequality and we obtain
\begin{eqnarray*}
C \frac{\G}{f (r)}\alpha'(r) & + & \G\alpha'' (r) + (n-1)\G\frac{f' (r)}{f
(r)}\alpha' (r)\\ & \ & \ \ -(n-1)\G\frac{1}{f^2 (r)} g(\alpha (r)) g' 
(\alpha
(r))> 0\
\end{eqnarray*}
for all $r \ge 1$. Since $g (y) > 0$ and
$g' (y) > 0$ for all $y \in {\real}^+$ and $\G > 0$, the above inequality
can be simplified  to
$$  {{\alpha'' (r)}\over {g(\alpha (r)) g' (\alpha (r))}}  + \frac{n-1}{f 
(r) }(C
 + f' (r))\, {{\alpha' (r) }\over { g(\alpha (r)) g' (\alpha (r))}}
-(n-1)\frac{1}{f^2 (r)} > 0\
\leqno (3.6)
$$  
for all $r \ge 1$. The proof proceeds similarly to \cite{Ch-Law-Leung} (cf.
\cite{R-R}). For
$r
\ge 1$, let 
$$ H (r) = {{\alpha' (r)}\over {g (\alpha (r)) g' (\alpha (r))}}\,.
$$  We have 
$$   {{\alpha'' (r)}\over {g(\alpha (r)) g' (\alpha (r))}} = H' (r) + H^2 (r)
[g'^2 (\alpha (r)) + g (\alpha (r)) g'' (\alpha (r))]
$$  for all $r \ge 1$. As $\alpha$ is bounded, we can find a constant 
$C_o > 1$
such that $C_o^2 \ge g'^2 (\alpha (r)) + g (\alpha (r)) g'' (\alpha (r))$ for
all $r \ge 1$. Then (3.6) becomes
$$  C_o^2 H^2 (r) + H' (r) +\frac{n-1}{f (r)}(C + f' (r))\, H (r)
-\frac{(n-1)}{f^2 (r)} > 0
\leqno (3.7)
$$ 
for all $r \ge 1$. For the quadratic form
$$  C_o^2 H^2 (r) + \frac{n-1}{f (r)}(C  + f' (r) ) H (r) -\frac{(n-1)}{f^2
(r)}\,,
$$  it is nonpositive for $H(r) \in [0,b_C(r)]$ where $b_C(r)$ is given by
$$  b_C (r)= {{n-1}\over {2C_o^2 f (r)}} \left\{ \left[ (C  +f' (r))^2 +
{{4C_o^2}\over {n-1}} \right]^{1/2} -(C + f' (r)) \right\} > 0\,. \leqno 
(3.8)
$$ 
Hence
$$ H'(r) > 0\quad\mbox{ whenever}\ \quad H(r)\leq b_C(r)\ .
\leqno (3.9)
$$  Consider the function 
$$ q (x) = \sqrt{x^2 + b^2} - x > 0\,,
$$  where $b$ is a positive number. If $x \le 0$, then $q (x) \ge b$. For 
$x >
0$, we have $q' (x) < 0$. Thus $q (x) \ge c_1^2$ if $|x|$ is bounded, where
$c_1$ is a positive number. As $|f'| \le b$, if $\alpha'$ is bounded, 
then $|C 
+ f' (r) |
\le C + b$ and
$f (r)
\le C_1 r$ for some positive constant $C_1 > 0$, we can find a positive 
number
$C_2$ such that   
$$
b_C (r) \ge {{C_2}\over r}
$$ 
for all $r \ge 1$. Thus 
$$\int_1^\infty b_C (r) dr = \infty\,.$$ On the other hand we have 
$$ - ( \ln {{g' (\alpha (r))}\over {g (\alpha (r))}} )' = H (r) [g'^2 (\alpha
(r)) - g (\alpha (r)) g'' (y)] \ge  c^2 H (r)
$$ 
for all $r \ge 1$. Thus 
$$ - \ln {{g' (\alpha (r))}\over {g (\alpha (r))}} \vert^\infty_1 = -
\int_1^\infty \left( \ln {{g' (\alpha (r))}\over {g (\alpha (r))}} 
\right)' dr
\ge c^2
\int_1^\infty H (r) dr\,.
$$ As $g' (r) > 0$ for all $r > 0$ and $y$ is bounded and positive for $r \ge
1$, we have 
$$
\int_1^{\infty} H(r) dr < \infty\ .
$$ Since $H \geq 0$ and $\int_1^{\infty} b_C (r)dr =\infty$ , for some
sufficiently large $r$ we have $H(r)\leq b_C(r)$ and $H'(r)\leq 0$ (cf.
\cite{R-R}), contradicting (3.9).\qed
Liouville's type theorems for harmonic functions are studied by Yau \cite{Y}.
Takakuwa \cite{Ta} studies Liouville's type theorems for $p$-harmonic 
maps and
Hong \cite{Hong} discusses Liouville's type theorems for exponentially 
harmonic
maps. Ratto and Rigoli \cite{R-R} study Liouville's type theorems for
rotationally symmetric harmonic maps. For $p$-harmonic maps, Liouville's type
theorems are studied in \cite{Ch-Law-Leung} and \cite{L2}.\\[0.2in]
{\bf Proposition 3.10.} \ \ {\it Assume that $g (r) = \sinh r$ and $0 \le 
f' (r)/
f(r) 
\le a$ for all $r \in {\real}^+$. Let $\alpha \in C^2 ((0\,, R))$ be a 
positive
solution to equation (1.1) with $\lim_{r \to 0^+} \alpha (r) = 0\,.$ If}
$${{d G}\over {dx}} (x) \ge c_o \,G (x) \ \ \ \ \ \ {\mbox{for \ \ all}} 
\ \ \ x
> 0\,,$$ 
{\it where $c_o$ is a positive constant, then there exists a positive 
constant
$C$ such that
$\alpha' (r)
\le C$ for all
$r
\in (R/2\,, R)$.}\\[0.1in] 
{\bf Proof.} \ \ For $n = 1$, the result follows from lemma 2.2. Assume 
that $n
\ge 2$, we have
\begin{eqnarray*}
&\,&{d \over {dr}} \G\\
 &=&{{d G}\over {dx}} \vert_{x = \theta (r)} {{d \theta
}\over {dr}} (r)\\
&=&{{d G}\over {dx}} \vert_{x = \theta (r)} \left\{ \alpha' (r)  \alpha''
(r) +  (n - 1) {{g (\alpha (r) ) g' (\alpha (r)) \alpha' (r)}
\over {f^2 (r)}} - (n - 1) {{g^2 (\alpha (r)) f' (r)}\over {f^3 (r)}}
\right\}.
\end{eqnarray*}
Substitute into equation (1.1) we have 
\begin{eqnarray*} 
(3.11) \ \ \ \ \ \ \ \ \ \ \  \ \ \ & \ & [\G + {{d G}\over {dx}} 
\vert_{x =
\theta (r)} (\alpha' (r))^2]
\alpha'' (r)\\
 =   & - & (n - 1) {{d G}\over {dx}} \vert_{x = \theta (r)} {{g (\alpha (r))
g' (\alpha (r))}\over {f^2 (r)}} (\alpha' (r))^2\\  
& - & (n
- 1) {{f'(r)}\over {f(r)}} \left[ \G - {{d G}\over {dx}} \vert_{x = 
\theta (r)}
{{g^2 (\alpha (r))}\over {f^2 (r)}} \right] \alpha' (r)\\
& + & (n - 1) \G {{ g
(\alpha (r)) g' (\alpha (r))}\over {f^2 (r)}}\,.\ \ \ \ \ \ \ \ \ \ \ \ \ 
\ \ \ 
\ \ \ \ \ \ \ \ \ \ \ \ \ \ \ \ \ \ \ \ 
\end{eqnarray*} 
As $0 \le f' (r)/f (r) \le a\,,$ $G' (x) \ge c_o G (x)$ and $g (y) g' (y) =
\sinh y \cosh y \ge (\sinh y)^2 = g^2 (y)\,,$ we have
\begin{eqnarray*}(3.12) \ \ \ \ \ \ \ \ \ \ \  & \ & [\G + {{d G}\over {dx}}
\vert_{x = \theta (r)} (\alpha' (r))^2]
\alpha'' (r)\\
& \le &- (n - 1) {{d G}\over {dx}} \vert_{x = \theta (r)} {{g (\alpha 
(r)) g'
(\alpha (r))}\over {f^2 (r)}} \left[ {1\over 2}(\alpha' (r))^2 - \alpha' (r)
\right]\\
&\ & \ \ \  - (n - 1) {{d G}\over {dx}} \vert_{x = \theta (r)} {{g 
(\alpha (r))
g' (\alpha (r))}\over {f^2 (r)}} \left[ {1\over 2} (\alpha' (r))^2 - {1\over
{c_o}}
\right]\,.\ \ \ \ \ \ \ \ \ \ \ \ \ \ \ \
\ \ \ \ \ 
\end{eqnarray*}
By lemma 2.1 we have $\alpha' (r) > 0$ for all $r \in (0, R)$. If at a point
$r \ge R/2$ we have
$$\alpha' (r) \ge \max \ \left\{\,2\,, \sqrt{2\over {c_o}}\, 
\right\}\,,$$ 
then (3.12) implies that
$$ 
[\G + {{d G}\over {dx}} \vert_{x = \theta (r)} (\alpha (r))^2] \alpha'' (r)
\le 0\,.
$$ 
That is, $\alpha'' (r) \le 0$ whenever $\alpha' (r) \ge  \max \ \{\,2\,,
\sqrt{2/{c_o}}\,\}$ and $r \ge R/2\,.$ Therefore we can find a positive 
constant
$C$ such that $\alpha' (r) \le C$ for all $r \ge R/2$.\qed It follows 
from the
above proposition that any local positive $C^2$-solution
$\alpha$ of equation (1.1) with $f$ and $g$ satisfy the conditions in 
proposition
3.10 (and 
$\lim_{r \to 0^+} \alpha (r) = 0$) can be extended to whole ${\real}^+$ (cf.
section 4). The same argument shows that in the above proposition, we can let
$R = \infty$ and there is a positive constant $C$ such that $\alpha' (r)
\le C$ for all $r \ge 1$. This contrasts with  rotationally symmetric
$p$-harmonic maps, where under the condition that $0 \le f' (r) \le 
c_o^2$ and
$f (r) \le c_1$ for all $r$ large and $p > 2$, a rotationally symmetric
$p$-harmonic map from such a function $f$ to the hyperbolic space with 
bounded
derivative has to be the constant map \cite{L2} (cf. \cite{R-R},
theorem 2.14).\\[0.2in] 
{\bf Theorem 3.13.} \ \ {\it Assume that there exist
positive constants $c_o, c_1\,, C_1\,, C_2$ and $r_o$ such that $c_o \le 
f (r)
\le c_1$,
$0
\le f' (r) \le C_1$ for all
$r \ge r_o$ and $C_2 < g' (y)$ for all $y > 0\,.$ Suppose that}
$$
{{d G}\over {dx}} (x) \le C_o G (x) \ \ \ \ \ \ {\mbox{for \ \ all}} \ \ 
\ x
> 0\,,
$$  
{\it where $C_o$ is a positive constant. If
$\alpha
\in C^2 ((0,
\infty))$ is a nonnegative solution to equation (1.1) with $f$ and $g$ 
satisfy
the above conditions and $\lim_{r \to 0^+} \alpha (r) = 0$ and $\lim_{r \to
\infty} \inf \alpha' (r) = 0$, then $\alpha \equiv 0$ on $(0, \infty)\,.$ In
particular, if
$\alpha$ is bounded, then $\alpha \equiv 0\,.$}\\[0.1in] 
{\bf Proof.}
\
\ Suppose that
$\alpha
\not\equiv 0$. Lemma 2.1 and the proof of lemma 2.7 imply that
$\alpha (r) > 0$ and $\alpha' (r) > 0$ for all $r > 0$. Hence there are 
positive
constants $C_3$ and $C_4$ such that $g (\alpha (r))
\ge C_3$ and $g' (\alpha (r)) \ge C_4$ for all $r \ge r_o$. For $r \ge 
r_o\,,$
as in (3.12) we have
\begin{eqnarray*}  
& \ & [\G + {{d G}\over {dx}} \vert_{x = \theta (r)} (\alpha'
(r))^2] \alpha'' (r)\\
& \ge & - (n - 1) C_o \G {{g (\alpha (r)) g' (\alpha (r))}\over {f^2 (r)}}
(\alpha' (r))^2 - (n - 1) {{f'(r)}\over {f(r)}} \G \alpha' (r)\\
& \ & \ \ \ + (n - 1) {{f'(r)}\over {f(r)}} {{d G}\over {dx}} \vert_{x = 
\theta
(r)} {{g^2 (\alpha (r))}\over {f^2 (r)}} \alpha' (r) +  (n - 1) \G {{
g (\alpha (r)) g' (\alpha (r))}\over {f^2 (r)}}\\
& \ge & (n - 1) \G \left\{  {{ g (\alpha (r)) g' (\alpha (r))}\over {f^2 
(r)}} [
1 -  C_o (\alpha' (r))^2 ] - {{C_1}\over {c_o}} \alpha' (r) \right\}\,.
\end{eqnarray*}
For $r \ge r_o$ we have
$$
{{ g (\alpha (r)) g' (\alpha (r))}\over {f^2 (r)}} \ge {{C_3 C_4}\over
{c_1^2}}\,.
$$ 
For $r \ge r_o$, if 
$$\alpha' (r) \le \min \ \left\{ {1\over {2 \sqrt{C_o}}}\,, \ {{c_oC_3 
C_4  
}\over {2c_1^2C_1}}
\right\}\,,$$ then
$$[\G + {{d G}\over {dx}} \vert_{x = \theta (r)} (\alpha' (r))^2] 
\alpha'' (r)
\ge 0\,$$
that is, $\alpha'' (r) \ge 0\,.$ Therefore we can find a positive constant
$\epsilon > 0$ such that $\alpha' (r) > \epsilon$ for all $r \ge r_o$. This
contradicts $\lim_{r \to \infty} \inf \alpha' (r) = 0\,.$ We note that if
$\alpha$ is bounded, then $\lim_{r \to \infty} \inf \alpha' (r) = 
0\,.${{\hfill
{$\rlap{$\sqcap$}\sqcup$}}} 
\vspace{0.5in}

{\bf \Large 4. \ \ Existence of Bounded Positive Solutions}

\vspace{0.3in}

As $G (x) > 0$ and $dG/dx (x) \ge 0$ for all $x \ge 0$, from (2.6) we have
\begin{eqnarray*}  (4.1) \ \ \ \ \ \ \ & \ & \theta' (r) \le  \theta' (r) 
{{\G + {{d G}\over {dx}}|_{x =
\theta (r)} (\alpha' (r))^2} \over {\G}} \\ & = & (n - 1) \{ {{2 g 
(\alpha (r)
) g' (\alpha (r))
\alpha' (r)}\over {f^2 (r)}} - {{f' (r)}\over {f (r)}} [ {{ g^2 (\alpha
(r))}\over {f^2 (r)}} + (\alpha' (r))^2 ] \}\,.\\  \ \ \ \ \ \ \  
\end{eqnarray*}
{\bf Lemma 4.2.} \ \ {\it Let $\alpha (r) \in C^2 (0, R)$ be a positive 
solution
to equation (1.1) with
$\lim_{r
\to 0^+} \alpha (r) = 0$, where
$R$ is a positive number. Suppose that there exists a positive constant 
$a$ such
that $0 < g' (y) \le a$ for all 
$y > 0$. Then
$\alpha$ can be extended to a positive solution of (1.1) on
$(0, \infty)$.}\\[0.1in]  {\bf Proof.} \ \ By lemma 2.2 we need only to
consider the case $n \ge 2$. From (4.1) we have
\begin{eqnarray*}
(4.3) \ \ \ \ \ \ \theta' (r) & \le & (n - 1) \left\{ {a\over {f (r)}} 
{{2 g
(\alpha (r))
\alpha' (r)}\over {f (r)}} - {{f' (r)}\over {f (r)}} [ {{ g^2 (\alpha 
(r))}\over
{f^2 (r)}} + (\alpha' (r))^2 ] \right\} \ \ \ \ \ \ \ \\ & \le & {C\over {f
(r)}} (a + |f' (r)|) \theta\,,
\end{eqnarray*} 
where $C$ is a positive constant. There exists a positive constant $b$  such
that 
$|f' (r)| \le b$ for all $r \in (0, R)$. An integration of (4.3) shows that
$\theta (r) < C$ for  all
$r \in (R/2, R)$. Thus $\alpha' (r)$ and hence $\alpha (r)$ are bounded on
$(R/2, R)$.  We can continue the solution to $(0, R + \epsilon)$ for some
positive  number
$\epsilon$ (\cite{CL}, p. 15). The proof is completed by a continuation
argument.{\hfill {$\rlap{$\sqcap$}\sqcup$}}\\[0.2in]
{\bf Theorem 4.4.} \ \ {\it For $n \ge 2$, assume that $\lim_{r \to 
\infty} f'
(r) =
\infty$ and
$0 < g' (y) \le a$ for all $y \ge 0$, where $a$ is a positive constant. Let
$\alpha \in C^2 ((0, \infty))$ be a positive solution to equation (1.1) with
$f$ and $g$ as above and $\lim_{r \to 0^+} \alpha (r) = 0$. Then for any
$\epsilon \in (0, \epsilon)$, there exist positive constants $r_o$ and 
$C$ such
that}
$$\alpha' (r) \le {{C}\over {f^{1 - \epsilon} (r)}}$$ {\it for all $r \ge 
r_o$.}
\\[0.1in] 
{\bf Proof.} \ \ Given $\epsilon \in (0, 1)$, we can find a positive
number
$r_o$ such that $f' (r) \ge a^2/\epsilon^2$ for all $r \ge r_o$. Then for
$r \ge r_o$, we have
\begin{eqnarray*} (4.5) \ \ \ & \ & {{2 g (\alpha (r)) g' (\alpha (r))
\alpha' (r)}\over {f^2 (r)}} - {{f' (r)}\over {f (r)}} [ {{ g^2 (\alpha
(r))}\over {f^2 (r)}} + (\alpha' (r))^2 ]\\ & = & {{2 g (\alpha (r)) g' 
(\alpha
(r))}\over {f^{3\over 2} (r)f'^{1\over 2} (r)}} [ {{f'(r) (\alpha' 
(r))^2}\over
{f (r)}} ]^{1\over 2} - {{ f' (r) g^2 (\alpha (r))}\over {f^3 (r)}} - {{f'
(r)}\over {f (r)}} (\alpha' (r))^2\\   & \le & {1\over {\epsilon}} {{g^2 
(\alpha
(r)) g'^2 (\alpha (r))}\over {f^3 (r) f' (r)}} + \epsilon {{f' (r)}\over {f
(r)}} (\alpha' (r))^2 - {{ f' (r) g^2 (\alpha (r))}\over {f^3 (r)}}
 - \,{{f' (r)}\over {f (r)}} (\alpha' (r))^2\\   & = & {{g^2
(\alpha (r)) }\over {\epsilon f^3 (r) f' (r)}} [g'^2 (\alpha (r)) - 
 \epsilon^2 f' (r)^2]\\  & \ & \ \ \ \ \ \ \ \ \  - \, (1 - \epsilon) [{{ 
f' (r)
g^2 (\alpha (r))}\over {f^3 (r)}} + {{f' (r)}\over {f (r)}} (\alpha' 
(r))^2]\\ 
& \le & {{ g^2 (\alpha (r))}\over {\epsilon f^3 (r) f' (r)}} [a^2 -
 \epsilon^2 f' (r)^2] - (1 - \epsilon) [{{ f' (r) g^2 (\alpha (r))}\over {f^3
(r)}} + {{f' (r)}\over {f (r)}} (\alpha' (r))^2] \ \ \  \\ &
\le & - (1 - \epsilon) {{f' (r)}\over {f (r)}} [{{ \alpha^2 (r)}\over 
{f^2 (r)}}
+ (\alpha' (r))^2]\,,
\end{eqnarray*}  
where we have used the inequality $2AB \le (1/\epsilon)A^2 + \epsilon 
B^2$ for
$\epsilon \in (0, 1)\,.$ By (4.1) and (1.2) we have
$$
\theta' (r) \le - 2 (1 - \epsilon) {{f' (r)}\over {f (r)}} \theta
$$ 
for all $r
\ge r_o$. An integration gives 
$$\alpha' (r) \le C f^{ - (1 - \epsilon)} (r)\,,$$ where $C$ is a positive
constant.{{\hfill {$\rlap{$\sqcap$}\sqcup$}}\\[0.2in]}
{\bf Corollary 4.6.} \ \ {\it For $n \ge 2$, assume that $\lim_{r \to 
\infty} f'
(r) =
\infty$ and
$0 < g' (y) \le a$ for all $y \ge 0$, where $a$ is a positive constant. 
If there
exist positive numbers $s > 1\,,$ $C'$ and $r'$ such that $f (r) \ge C' 
r^s$ for
all $r \ge r'$, then any positive $C^2$-solution $\alpha$ of equation 
(1.1) with
$\lim_{r \to 0^+} \alpha (r) = 0$ is bounded.}\\[0.1in] {\bf Proof.} \ \ By
lemma (4.2) we may assume that $\alpha$ is defined on $(0,
\infty)$. Choose $\epsilon \in (0, 1)$ such that $(1 - \epsilon) s > 1$. 
Then we
can find a positive number $r_o$ such that
$$
\alpha' (r) \le {C \over {f^{ (1 - \epsilon)} (r)}} \le {{C''}\over {r^{ 
(1 -
\epsilon)s} }} \leqno (4.7)
$$ for all $r \ge r_o$, where $C''$ is a positive constant. By 
integrating (4.7)
we conclude that
$\alpha$ is bounded.\qed 
In particular, any rotationally symmetric $F$-harmonic maps from the
hyperbolic space to the Euclidean space is bounded. We consider rotationally
symmetric $F$-harmonic maps from the hyperbolic space to
itself.\\[0.2in] {\bf Lemma 4.8.} \ \ {\it Assume that in equation (1.1) $f
\equiv g$  and $f'' (r) = g'' (r) \ge 0$ on
${\real}^+$. Let $\alpha (r) \in C^2 (0, R)$ be a positive solution to 
equation
(1.1) with
$\lim_{r
\to 0^+} \alpha (r) = 0$, where
$R$ is a positive number. If there is a point $r_o \in (0, R)$ such that
$\alpha (r_o) \le r_o$, then $\alpha$ can be extended to a positive 
solution of
(1.4) on
$(0, \infty)$.}\\[0.1in]   
{\bf Proof.} \ \ We may assume that $\alpha
\not\equiv 0$ on $(0, R)$. By proposition 2.24 we have
$\alpha (r)
\le r$ for all
$r
\in (0, R)$. Thus $\alpha (r) \le R$ for all $r \in (0, R)$. Since 
$\alpha' (r) >
0$, $\lim_{r \to R^-} \alpha (r)$ exists.  As in lemma 4.2, $\alpha$ can be
extended to a positive solution on $(0, \infty)$.{{\hfill
{$\rlap{$\sqcap$}\sqcup$}}\\[0.2in]} 
{\bf Lemma 4.9.} \ \ {\it Let $\alpha \in C^2 (0, \infty)$ be a positive 
solution
to equation (1.1) with $f (r) = g (r) =
\sinh r$ on $(0, \infty)$ and $\lim_{r \to 0^+} \alpha (r) = 0$. If there 
exists
a positive number $r_o > (\ln 3)/2$ such that $\alpha (r_o) = r_o -
\delta$ for some $\delta \ge 0$ and $\alpha' (r_o) < 1$, then $\alpha (r) 
< r_o
- \delta$ and $\alpha' (r) < 1$ for all $r > r_o$.}\\[0.1in] 
{\bf Proof.} \ \ For $n = 1$, the result follows from lemma 2.2. We may 
assume
that $n \ge 2$.  Suppose that there is a point $r' > r_o$ such that 
$\alpha (r')
= r' -
\delta$. We can find  a point ${\tilde r} \in (r_o,
 r')$ such that $\alpha ({\tilde r}) < {\tilde r} - \delta$, $\alpha' 
({\tilde
r}) = 1$ and $\alpha'' ({\tilde r})
\ge 0$. Let $\alpha ({\tilde r}) = {\tilde r} - \epsilon\,,$ where 
$\epsilon >
\delta$ is a positive constant. As $f ({\tilde r}) = g ({\tilde r}) = \sinh
\,{\tilde r}$ and at the point ${\tilde r}$ where
$\alpha' ({\tilde r}) = 1\,,$ we have
\begin{eqnarray*}(4.10) & \ & \ {{2 g (\alpha ({\tilde r})) g' (\alpha
({\tilde r}))
\alpha' ({\tilde r})}\over {f^2 ({\tilde r})}}  - {{f' ({\tilde r})}\over {f
({\tilde r})}} [ {{g^2 (\alpha ({\tilde r}))}\over {f^2 ({\tilde r})}} +
(\alpha' ({\tilde r}))^2 ]\\ & = & {1\over 8 f^3 ({\tilde r})} [ 2 
(e^{\tilde r}
- e^{-{\tilde r}}) (e^\alpha - e^{-\alpha}) (e^\alpha + e^{-\alpha}) \\ & 
\ & \
\ \ \ \ \ \ \ \ \ \ \ \  - (e^{\tilde r} + e^{-{\tilde r}}) (e^\alpha -
e^{-\alpha})^2 - (e^{\tilde r} + e^{-{\tilde r}}) (e^{\tilde r} - e^{-{\tilde
r}})^2 ]\\  & = & {1\over {8 f^3 ({\tilde r})}} (e^{\tilde r} e^{2\alpha} -
3e^{\tilde r} e^{- 2\alpha} - 3 e^{-{\tilde r}} e^{2\alpha} + e^{-{\tilde r}}
e^{-2\alpha} + 3 e^{\tilde r} + 3 e^{-{\tilde r}} - e^{-3{\tilde r}} -
e^{3{\tilde r}})\\ & = & {1\over 8 f^3} [ (e^{-2\epsilon} - 1) ( 
e^{3{\tilde r}}
- 3 e^{\tilde r}) + (3 e^{-{\tilde r}} - e^{-3{\tilde r}}) (1 - 
e^{2\epsilon})
]\,, \ \ \ \ \ \ \  \ \ \ \ \ \ \ 
\end{eqnarray*}  as $\alpha ({\tilde r}) = {\tilde r} - \epsilon$.  We 
have $3
e^{-{\tilde r}} - e^{-3{\tilde r}} > 0\,.$ Since ${\tilde r} > {\tilde 
r}_o >
(\ln 3)/2$, we have $e^{3{\tilde r}} - 3 e^{\tilde r} > 0$  as well. 
Therefore
(4.10) implies that 
$$ {{2 g (\alpha ({\tilde r})) g' (\alpha ({\tilde r})) \alpha' ({\tilde
r})}\over {f^2 ({\tilde r})}}  - {{f' ({\tilde r})}\over {f ({\tilde 
r})}} [
{{g^2 (\alpha ({\tilde r}))}\over {f^2 ({\tilde r})}} + (\alpha' ({\tilde 
r}))^2
] < 0\,.
$$  
By (4.1) we have $\theta' ({\tilde r}) < 0$. On the other hand we have
\begin{eqnarray*}
\theta' ({\tilde r}) &  = & (n - 1) ( {{g^2 (\alpha ({\tilde r}))} \over {f^2
({\tilde r})}} )' + 2 \alpha' ({\tilde r})
\alpha'' ({\tilde r})\\ & = & (n - 1) {{2 f ({\tilde r}) g (\alpha 
({\tilde r}))
[ f({\tilde r}) g' (\alpha) \alpha' ({\tilde r}) - f' ({\tilde r}) g (\alpha)
]}\over {f^4 ({\tilde r})}} + 2 \alpha'' ({\tilde r})\\ & = & (n - 1) {{g
(\alpha)}\over {f^3 ({\tilde r})}} (e^{\tilde r} e^{-\alpha ({\tilde r})} -
e^{-{\tilde r}} e^{\alpha ({\tilde r})} ) + 2
\alpha'' ({\tilde r})\\ & > & 2 \alpha'' ({\tilde r})\,,
\end{eqnarray*} 
as $\alpha ({\tilde r}) < {\tilde r}$ and $\alpha' ({\tilde r})
= 1\,.$ Thus $\alpha'' ({\tilde r}) < 0\,,$ which is a contradiction. 
Therefore
$\alpha (r) < r - \delta$ for all $r > r_o$. Similarly, we can show that
$\alpha' (r) < 1$ for all
$r >r_o$.{{\hfill{$\rlap{$\sqcap$}\sqcup$}}\\[0.2in]}     
{\bf Lemma 4.11.} \ \
{\it For $n \ge 2$, let $\alpha \in C^2 (0, \infty)$ be a positive 
solution to
equation (1.1) with $f (r) = g (r) =
\sinh r$ on $(0, \infty)$ and $\lim_{r \to 0^+} \alpha (r) = 0$. If there 
exists
a positive number $r_o > (\ln 3)/2$ such that $\alpha (r_o) < r_o$ and
$\alpha' (r_o) < 1$, then $\alpha$ is bounded on $(0, \infty)$.}\\[0.1in] 
{\bf Proof.}
\
\ Let $\alpha (r_o) = r_o - \delta$. Then by lemma 4.9, we have $\alpha (r)
< r - \delta$ for all $r > r_o$. As in (4.5) (with
$\tau = 1 -
\epsilon$), we have
\begin{eqnarray*}(4.12) \ \ \ \ \ \ \ \ \ \ \  & \ & \ \ {{2 g (\alpha 
({\tilde
r})) g' (\alpha ({\tilde r}))
\alpha' ({\tilde r})}\over {f^2 ({\tilde r})}}  - {{f' ({\tilde r})}\over {f
({\tilde r})}} [ {{g^2 (\alpha ({\tilde r}))}\over {f^2 ({\tilde r})}} +
(\alpha' ({\tilde r}))^2 ]\\  & = & {{g^2 (\alpha (r)) }\over {(1 - \tau) 
f^3 
(r) f' (r)}} [g'^2 (\alpha (r))
- 
 (1 - \tau)^2 f' (r)^2]\\  & \ & \ \ \ \ \ \ \ \ \  - \, \tau  [{{ f' (r)
g^2 (\alpha (r))}\over {f^3 (r)}} + {{f' (r)}\over {f (r)}} (\alpha'
(r))^2]\,.\ \ \ \ \ \ \ \ \ \ \ \ \ \ \ \ \ \ \ \ \ \ \ \ \ \ \ \ \ \ \ 
\end{eqnarray*}
If we choose $\tau \in (0, 1)$ such that $(1 - \tau)^2 \ge (1
+ e^{-2
\delta})/2$, then a calculation shows that (cf. \cite{L1})
$$
{{g^2 (\alpha
(r)) }\over {(1 - \tau) f^3 (r) f' (r)}} [g'^2 (\alpha (r)) - 
 (1 - \tau)^2 f' (r)^2] < 0\
$$
for all $r > r_o$. As $f' (r) = \cosh r\,,$ we have $c_o < f' (r)/ f (r)$ for
some positive constant $c_o\,.$ By (4.1) and (4.12) we can find a positive
constant
$c$ such that 
$$
\theta' (r) \le - c
\theta (r)
$$
for all $r > r_o$  Therefore $\theta (r)$ and hence $\alpha' (r)$ decay
exponentially. Hence
$\alpha$ is bounded.{\mbox{\bf \ \ \ \ \ Q.E.D.}\\[0.2in]}
{\bf Lemma 4.13.} \ \  \ \ {\it For $n \ge 2$, let $\alpha \in C^2 (0, 
\infty)$
be a positive solution to equation (1.1) with $f (r) = g (r) =
\sinh r$ on $(0, \infty)$ and 
$\lim_{r \to 0^+} \alpha (r) = 0$. If there exists a point $r_o > 0$  
such that
$\alpha (r_o) \le r_o$, then either $\alpha$ is bounded or $\alpha$ is
asymptotic to the identity map $\alpha_I (r) = r\,.$}\\[0.1in]
{\bf Proof.} \ \ By proposition 2.24, we have $\alpha (r) = r$ or $\alpha (r)
< r$ for all $r > 0$. Assume that $\alpha (r) < r$ for all $r > 0$. If 
$\alpha'
(r) < 1$ for some $r > (\ln 3)/2$, then lemma 4.11 implies that
$\alpha$ is bounded on $(0, \infty)$. Assume that $\alpha' (r) \ge 1$ for all
$r >  (\ln 3)/2$. Then the function
$$
\rho (r) = r - \alpha (r)
$$ 
is nonincreasing for $r > (\ln 3)/2$. Therefore $\lim_{r \to \infty} \rho (r)
= c$, where $c$ is a nonnegative number. If $c > 0$, then we can find a
positive number $r'$ such that $\rho (r) = r - \alpha (r) > c/2$ for all 
$r >
r'$, that is, $\alpha (r) < r - c/2$. The proof of lemma 4.11 shows that
$\alpha$ is a bounded function. If $c = 0$, then $\alpha$ is asymptotic 
to the
identity map $\alpha_I (r) = r\,.${{\hfill 
{$\rlap{$\sqcap$}\sqcup$}}\\[0.2in]}
{\bf Theorem 4.14.} \ \  \ \ {\it For $n \ge 2$, let $\alpha \in C^2 (0,
\infty)$ be a positive solution to equation (1.1) with $f (r) = g (r) = 
\sinh r$
and 
$\lim_{r \to 0^+} \alpha (r) = 0\,.$ If there exists a point $r_o > 0$  
such that
$\alpha (r_o) < r_o$, then $\alpha$ is bounded.}\\[0.1in]
{\bf Proof.} \ \ Suppose that $\alpha$ is asymptotic to
the identity map
$\alpha_I (r) = r\,.$ As $\alpha_I (r) > \alpha (r)$ for all $r > 0$ and 
$$
\lim_{r\to 0^+} \left[ \alpha_I (r) - \alpha (r) \right] = 0  \ \ \ \
{\mbox{and}} \ \ \ 
\lim_{r\to \infty} \left[ \alpha_I (r) - \alpha (r) \right] = 0\,,
$$
there exists a point $\bar r > 0$ such that 
$$
\alpha_I (\bar r) > \alpha (\bar r)\,, \ \ \ \alpha'_I (\bar r) = \alpha' 
(\bar
r) = 1\,, \ \ \ {\mbox{and}} \ \ \ \alpha''_I (\bar r) \le \alpha'' (\bar
r)\,.
$$
That is,
$$
\bar r > \alpha (\bar r)\,, \ \ \ \alpha' (\bar r) = 1\,,
\ \ \ {\mbox{and}} \ \ \  \alpha'' (\bar r) \ge 0\,.\leqno (4.15)
$$
From (3.11) we have
\begin{eqnarray*}
(4.16) \ & \ & [\G + {{d G}\over {dx}} \vert_{x =
\theta (r)} (\alpha' (r))^2]
\alpha'' (r)\\
 & =  & \ \,\,- (n - 1) {{d G}\over {dx}} \vert_{x = \theta (r)} \left[ {{g
(\alpha (r)) g' (\alpha (r))}\over {f^2 (r)}} (\alpha' (r))^2 - {{f' 
(r)}\over {f
(r)}} {{g^2 (\alpha (r))}\over {f^2 (r)}} \alpha' (\bar r) \right]\\
& \  & \ \ - (n - 1) \G \left[ {{f'(r)}\over {f(r)}} \alpha' (\bar r) - 
{{ g
(\alpha (r)) g' (\alpha (r))}\over {f^2 (r)}} \right]\\
& = & \ \,\,- (n - 1) {{d G}\over {dx}} \vert_{x = \theta (r)} {{\sinh 
\,(\alpha
(\bar r)) }\over {\sinh^3 \bar r}} \left[ \cosh \,(\alpha (\bar r)) \sinh 
\bar r 
- \sinh \,(\alpha (\bar r)) \cosh \bar r \right]\\
& \ & \ \ - (n - 1) \G \left[ {{\cosh \bar r}\over {\sinh
\bar r}} - {{\sinh \,(\alpha (\bar r)) \cosh
\,(\alpha (\bar r)) }\over {\sinh^2 \bar r}} \right]\\
& < & 0\,,
\end{eqnarray*}
as $\cosh \,(\alpha (\bar r)) \sinh \bar r  - \sinh \,\alpha (\bar r) 
\cosh \bar
r = [e^{-\alpha (\bar r)} e^{\bar r} - e^{\alpha (\bar r)} e^{- \bar 
r}]/2 > 0$
because $\alpha (\bar r) < \bar r\,.$ From (4.16) we have $\alpha'' (\bar 
r) <
0$, which contradicts (4.15). Therefore $\alpha$ cannot be asymptotic
to the identity map
$\alpha_I (r) = r\,.$ By lemma 4.13, $\alpha$ is a bounded function on $(0,
\infty)$.\qed

\pagebreak

\end{document}